\documentclass[11pt,reqno]{amsart}
\usepackage[T1]{fontenc}
\usepackage{lmodern}
\usepackage{amsmath,amssymb,mathtools,mathrsfs}
\usepackage[margin=28mm]{geometry}
\usepackage{tikz}
\usepackage{pgfplots}
\usepgfplotslibrary{groupplots}
\pgfplotsset{compat=1.18}
\usepackage{xcolor}
\usepackage{placeins}
\usepackage[numbers,sort&compress]{natbib}
\usepackage[hidelinks]{hyperref}
\allowdisplaybreaks[2]
\numberwithin{equation}{section}
\newtheorem{theorem}{Theorem}[section]
\newtheorem{proposition}[theorem]{Proposition}
\newtheorem{lemma}[theorem]{Lemma}
\newtheorem{corollary}[theorem]{Corollary}
\theoremstyle{remark}

\newcommand{\R}{\mathbb R}
\newcommand{\E}{\mathbb E}
\newcommand{\ind}{\mathbf 1}
\newcommand{\norm}[1]{\left\|#1\right\|}
\title[Local linearization and parameter estimation]
{Bi-parameter local linearization and parameter estimation
for the damped stochastic Klein--Gordon equation
with spatially colored noise}

\author[Y. Li]{Yumeng Li}
\address{Yumeng Li, School of Statistics and Mathematics,
Zhongnan University of Economics and Law,
Wuhan, Hubei 430073,   China}
\email{li-yu-meng@163.com}

\author[R. Wang]{Ran Wang}
\thanks{Corresponding author: Ran Wang (\texttt{rwang@whu.edu.cn}).}
\address{Ran Wang, School of Mathematics and Statistics, Wuhan University, Wuhan 430072, China}
\email{rwang@whu.edu.cn}
\date{}
\hypersetup{pdftitle={Bi-parameter local linearization and parameter estimation for the damped stochastic Klein-Gordon equation with spatially colored noise},pdfauthor={Yumeng Li, Ran Wang}}
\begin{document}
\begin{abstract}
Consider the one-dimensional nonlinear damped stochastic Klein--Gordon
equation driven by Gaussian noise that is white in time and fractional
in space, with Hurst parameter $H\in(1/2,1)$. For a globally Lipschitz
diffusion coefficient, we approximate characteristic mixed increments
by those of the additive equation, multiplied by the coefficient frozen
at the base point. A critical-kernel reduction and localized Riesz
covariance estimates control the additional kernel terms and spatial
correlations. We also establish $L^p$ error bounds for normalized
quadratic variations and, under a nondegeneracy condition, construct
jointly strongly consistent estimators of the Hurst parameter and a
positive diffusion parameter.
\end{abstract}

\maketitle

\section{Introduction and main results}

Consider the nonlinear damped stochastic Klein--Gordon equation
\begin{equation}\label{eq:model}
\begin{cases}
\bigl(\partial_t^2-\partial_x^2+a\partial_t+m^2\bigr)u(t,x)
 =F(u(t,x))\dot W(t,x),&t>0,\ x\in\mathbb R,\\
u(0,x)=\partial_tu(0,x)=0,&x\in\mathbb R.
\end{cases}
\end{equation}
Here, $a\in\R$ and $m\ge0$ are the damping and mass parameters,
and $F:\R\to\R$ is globally Lipschitz.
The driving noise $\dot W$ is a centered generalized Gaussian
field on a probability space $(\Omega,\mathcal F,\mathbb P)$,
white in time and spatially homogeneous with fractional spatial
covariance.
More precisely, it is specified by a centered Gaussian family
$\{W(\varphi):\varphi\in\mathcal D\}$, indexed by the space
$\mathcal D:=C_0^\infty([0,\infty)\times\mathbb R)$ of smooth,
compactly supported test functions, with covariance
\begin{equation}\label{Eq:noise_cov}
\E[W(\varphi)W(\psi)]
=
\int_0^\infty
\int_{\mathbb R}
\mathscr{F}(\varphi(t,\cdot))(\xi)\,
\overline{\mathscr{F}(\psi(t,\cdot))(\xi)}
\,\mu(d\xi)\,dt,
\qquad \varphi,\psi\in\mathcal D.
\end{equation}
The spatial spectral measure is
\begin{equation}\label{eq:spectral-density}
\mu(d\xi)=c_H|\xi|^{1-2H}\,d\xi,
\qquad
c_H=\frac{\Gamma(2H+1)\sin(\pi H)}{2\pi},
\qquad H\in (1/2,1),
\end{equation}
and $\mathscr{F}$ denotes the Fourier transform in the spatial
variable, with the convention
\[
\mathscr{F}(\varphi(t,\cdot))(\xi)
:=
\int_{\mathbb R}e^{-i\xi y}\varphi(t,y)\,dy.
\]
The spectral density in \eqref{eq:spectral-density} reduces to
$(2\pi)^{-1}$ at $H=1/2$, corresponding to space-time white noise. When $a=m=0$,
\eqref{eq:model} reduces to the stochastic wave equation.

The stochastic integration theory for spatially homogeneous noise
was developed by Dalang \cite{Dalang1999}, extending the
martingale-measure framework of Walsh \cite{Walsh1986}; see also
\cite{DKMNX09,DM2003,K14}.

\subsection{Local linearization and related work}
Local linearization describes the small-scale fluctuations of
many parabolic stochastic partial differential equations (SPDEs)
with multiplicative noise; see, for example,
\cite{FKM2015,HP15}. Let $X$ solve the one-dimensional stochastic
heat equation with globally Lipschitz coefficient $F$, zero initial
data, and space-time white noise, and let $Z$ be its additive
counterpart driven by the same noise. For fixed $t>0$ and $x\in\R$,
Foondun et al.~\cite{FKM2015} showed that
\begin{equation}\label{SHE appro}
\begin{aligned}
X(t,x+h)-X(t,x)
={}&F(X(t,x))\bigl[Z(t,x+h)-Z(t,x)\bigr]
   +\mathcal R_{t,x}(h),
\end{aligned}
\end{equation}
where
\[
 \|\mathcal R_{t,x}(h)\|_{L^2(\Omega)}
 =o(|h|^{1/2}),\qquad h\to0.
\]
Thus the Gaussian increment governs the leading fluctuation,
with the diffusion coefficient frozen at $(t,x)$.
This approximation is related to Gubinelli's controlled-path
viewpoint \cite{G04}. Temporal analogues were studied for space-time white noise in
\cite{KSXZ2013,HP15} and for spatially fractional noise in
\cite{WX2024,QWWX2025}; spatial approximations with Riesz noise
were obtained in \cite{W2024}.
 These local linearization results and related Gaussian
approximations are used to study laws of the iterated logarithm,
power variations, and small-ball probabilities for parabolic SPDEs
\cite{FKM2015,PT2007,WX2024,QWWX2025,KKM23}.
Related work concerns multifractal behavior \cite{HK2017}, dense blowup
\cite{CHKK19}, hitting probabilities \cite{DNP2025}, Hausdorff
dimensions \cite{KLPX25}, and local properties of the
Kardar--Parisi--Zhang equation
\cite{Das2024,HL2025}.

For stochastic wave equations, ordinary time or space increments
involve noise contributions from thin strips along the backward
light cone, which do not shrink to the observation point.
Huang and Khoshnevisan \cite{HK2016} showed, for a nonlinear wave
equation driven by space-time white noise, that the limiting
quadratic variations depend on the diffusion coefficient over
extended characteristic regions. Thus the parabolic local
linearization obtained by freezing the coefficient at the
observation point generally fails.
Huang, Oh, and Okamoto \cite{HOO2024} established local linearization
in the space-time white-noise case by using mixed second-order
increments along the two characteristic directions, which localize
the stochastic integral to a small region. Liu and Wang
\cite{LW2025} and Li and Wang \cite{LiWang2026} extended this
approach to spatially fractional noise in the regimes
$H\in(1/2,1)$ and $H\in(1/4,1/2)$, respectively.

Rang and Wang \cite{RW2026} obtained local linearization for the
nonlinear damped stochastic Klein--Gordon equation driven by
space-time white noise, using the critical-kernel reduction of
Chen and Lee \cite{CL2026}.  In this paper, we extend this analysis to spatially
colored noise with $H\in(1/2,1)$. The critical mixed kernel leaves
contributions on the old cone and characteristic strips.
Noise integrals over disjoint spatial regions can be correlated
under the Riesz covariance.
We control these contributions in the covariance norm and handle
the zeroth-order correction through a smoother drift convolution. This gives
the same leading Gaussian fluctuation scale as for the wave
equation, together with a quantitative remainder estimate.

\subsection{Characteristic increments and local linearization}

To formulate the local linearization result, let $U$ be the solution
of the additive damped stochastic Klein--Gordon equation
\begin{equation}\label{eq:linear-model}
\begin{cases}
(\partial_t^2-\partial_x^2+a\partial_t+m^2)U(t,x)
 =\dot W(t,x),&t>0,\ x\in\R,\\
U(0,x)=\partial_tU(0,x)=0,&x\in\R.
\end{cases}
\end{equation}
Here, $U$ and $u$ are driven by the same
noise $W$, with the same parameters $a$, $m$, and $H$.
Following \cite{HOO2024,LW2025,RW2026}, we introduce the
characteristic coordinates
\begin{equation}\label{eq transform}
(\tau,\lambda)
 =\left(\frac{t-x}{\sqrt2},\frac{t+x}{\sqrt2}\right).
\end{equation}
The two coordinate directions are aligned with the characteristics
of the wave operator. For $z=(\tau,\lambda)$ with
$\tau+\lambda\ge0$, write
\begin{equation}\label{eq:null-coordinates}
\begin{aligned}
t_z&=\frac{\tau+\lambda}{\sqrt2},&
x_z&=\frac{\lambda-\tau}{\sqrt2},\\
v(z)&=u(t_z,x_z),&V(z)&=U(t_z,x_z).
\end{aligned}
\end{equation}

For a signed step $r$, define
\begin{equation}\label{eq diff}
\begin{aligned}
\delta_r^{(1)}f(\tau,\lambda)
 &=f(\tau+r,\lambda)-f(\tau,\lambda),\\
\delta_r^{(2)}f(\tau,\lambda)
 &=f(\tau,\lambda+r)-f(\tau,\lambda),
\end{aligned}
\end{equation}
whenever the evaluation points belong to the half-plane
$\tau+\lambda\ge0$. For $h>0$, set
$D_h^\pm=\delta_{\pm h}^{(1)}\delta_h^{(2)}$ and abbreviate
$D_h=D_h^+$. In particular,
\begin{equation}\label{eq:rectangular-increment}
\begin{aligned}
D_hf(\tau,\lambda)
=f(\tau+h,\lambda+h)-f(\tau+h,\lambda)-f(\tau,\lambda+h)+f(\tau,\lambda).
\end{aligned}
\end{equation}
The corresponding local linearization remainders are
\begin{equation}\label{eq R}
 R_h^\pm(z)=D_h^\pm v(z)-F(v(z))D_h^\pm V(z).
\end{equation}
Fix $T>0$. We call $(z,h)$ admissible for the positive orientation
if $0\le t_z$ and $t_z+\sqrt2h\le T$; for the negative orientation,
we require $h/\sqrt2\le t_z\le T-h/\sqrt2$.
These conditions ensure that all four vertices correspond to times
in $[0,T]$.

\begin{theorem}[Local linearization]\label{thm:linearization}
For every $T>0$ and $1\le p<\infty$, there exists a finite constant
$C=C(p,T,a,m,H,F)$ such that
\begin{equation}\label{eq:main-linearization}
 \|R_h^\pm(z)\|_{L^p(\Omega)}\le C h^{H+1}
\end{equation}
for all admissible $(z,h)$ with $0<h\le1$. The constant is
uniform in $x_z\in\R$.
\end{theorem}

For the additive Klein--Gordon solution $V$, \eqref{eq:linear-cell}
and Lemma~\ref{lem:cell-variance} imply
\[
 \|D_h^\pm V(z)\|_2\sim\sqrt{\kappa_H}\,h^{H+1/2},
 \qquad h\downarrow0,
\]
where
\begin{equation}\label{eq:kappa-intro}
 \kappa_H=\frac{2^{H-3/2}}{2H+1}.
\end{equation}
The constant $\kappa_H$ is the characteristic-cell variance coefficient
computed by Liu and Wang~\cite[Lemma~3.1, equation~(3.1)]{LW2025}. Gaussianity and
\eqref{eq:main-linearization} therefore yield
\[
 \frac{\|R_h^\pm(z)\|_{L^p(\Omega)}}
      {\|D_h^\pm V(z)\|_{L^p(\Omega)}}
 =O(h^{1/2}),\qquad h\downarrow0.
\]
\subsection{Quadratic variations and statistical application}

We next study the quadratic variation of $v$.
Related results on variations and parameter estimation for SPDEs
can be found in \cite{KTZ2018,AGT2022,KT2018,SST2020,Shev2023,TKU2026,TZ2025}.

For each integer $N\ge1$, define
\[
 \tau_i=\frac{i}{N},\qquad \lambda_j=\frac{j}{N},\qquad
 z_{ij}=(\tau_i,\lambda_j),\qquad 0\le i,j\le N.
\]
On the observation square $[0,1]^2$, let
\begin{equation}\label{eq:grid-increment}
 \Delta_{ij}^N v:=D_{1/N}v(z_{ij}),\qquad 0\le i,j<N.
\end{equation}
The normalized quadratic variation is
\begin{equation}\label{eq:variation-def}
 Q_N(v)=N^{2H-1}
 \sum_{i=0}^{N-1}\sum_{j=0}^{N-1}\left(\Delta_{ij}^N v\right)^2.
\end{equation}

\begin{theorem}[Quadratic variations]\label{thm:variation}
For every $1\le p<\infty$, there exists a finite constant
$C_p=C(p,a,m,H,F)$ such that
\begin{equation}\label{eq:main-variation}
\left\|Q_N(v)-\kappa_H
 \int_0^1\int_0^1 F^2(v(\tau,\lambda))\,d\tau\,d\lambda
 \right\|_{L^p(\Omega)}
 \le C_pN^{-1/2},\qquad N\ge1.
\end{equation}
Moreover,
\begin{equation}\label{eq:variation-as}
Q_N(v)\longrightarrow\kappa_H
 \int_0^1\int_0^1 F^2(v(\tau,\lambda))\,d\tau\,d\lambda
 \qquad\text{almost surely as }N\to\infty.
\end{equation}
\end{theorem}

As an application, we construct jointly strongly consistent
estimators of the Hurst parameter and a positive diffusion
parameter, with error bounds in probability. Their construction
and proofs are given in
Sections~\ref{sec:estimation}--\ref{sec:hurst-estimation}.

\subsection{Organization}
We organize the remainder of the paper as follows.
In Section~\ref{sec:prelim}, we introduce the stochastic integration
framework and establish well-posedness and regularity.
In Section~\ref{sec:local}, we recall the critical-kernel reduction
and prove the local linearization theorem.
In Section~\ref{qv:section}, we establish the quadratic-variation
limits and apply them to estimation of the Hurst and diffusion
parameters. Finally, we illustrate their joint estimation
numerically in Section~\ref{sec:simulation}.

\section{Preliminaries and properties of the solution}\label{sec:prelim}
\subsection{The Riesz covariance and stochastic integration}
For $H\in(1/2,1)$, the covariance \eqref{Eq:noise_cov} can
equivalently be written in physical space as
\begin{equation}\label{eq:noise-intro}
 \E[W(\varphi)W(\psi)]
 =H(2H-1)\int_0^\infty\int_{\R^2}
 \varphi(s,y)\psi(s,y')|y-y'|^{2H-2}\,dy\,dy'\,ds.
\end{equation}
The spatial covariance defines a Hilbert space $\mathcal H$,
obtained by completing the real-valued
functions in $C_c^\infty(\R)$ under the inner product
\begin{equation}\label{eq:covariance-inner-product}
\begin{aligned}
 \langle\varphi,\psi\rangle_{\mathcal H}
 &=H(2H-1)\int_{\R^2}
   \varphi(y)\psi(y')|y-y'|^{2H-2}\,dy\,dy'\\
 &=\int_\R\widehat\varphi(\xi)
          \overline{\widehat\psi(\xi)}\,\mu(d\xi).
\end{aligned}
\end{equation}
Here $\widehat\varphi=\mathscr F\varphi$ uses the convention in
the introduction. The spectral measure \eqref{eq:spectral-density}
satisfies Dalang's condition
\[
 \int_\R\frac{\mu(d\xi)}{1+\xi^2}<\infty.
\]

The covariance \eqref{Eq:noise_cov} extends $W$ by isometry to
$L^2(\R_+;\mathcal H)$. Let $(\mathcal F_t)_{t\ge0}$ be its completed, right-continuous
natural filtration. For a predictable $\mathcal H$-valued process $g$
such that $\E\int_0^\infty\|g(s,\cdot)\|_{\mathcal H}^2\,ds<\infty$,
the stochastic integral satisfies
\begin{equation}\label{eq:stochastic-isometry}
 \E\left|\int_0^\infty\int_\R g(s,y)\,W(ds,dy)\right|^2
 =\E\int_0^\infty\|g(s,\cdot)\|_{\mathcal H}^2\,ds.
\end{equation}
See \cite{Walsh1986,Dalang1999}.

For $p\ge2$, a deterministic kernel $k$, and a predictable random
field $Z$, whenever the right-hand side below is finite, the
Burkholder--Davis--Gundy inequality gives
\begin{equation}\label{eq:bdg}
\begin{aligned}
 &\left\|\int_0^\infty\int_\R
          k(s,y)Z(s,y)\,W(ds,dy)\right\|_p^2\\
  \le &\, C_pH(2H-1)\int_0^\infty\int_{\R^2}
 |k(s,y)|\,|k(s,y')|\,|y-y'|^{2H-2}\|Z(s,y)\|_p\|Z(s,y')\|_p\,dy\,dy'\,ds.
\end{aligned}
\end{equation} 
Here $\|X\|_p=(\E|X|^p)^{1/p}$ for $1\le p<\infty$.
 
\subsection{Moment bounds and regularity}
Put $b=a^2/4-m^2$ and extend causal kernels by zero to $t\le0$.
The Green function of the original operator has the physical-space
representation
\begin{equation}\label{eq:fullgreen}
 G_{a,m}(t,x)=\frac12e^{-at/2}\ind_{\{|x|<t\}}\,
 \mathcal B_b\left(t^2-x^2\right),\qquad t>0,
\end{equation}
where
\[
 \mathcal B_b(r)=\sum_{k=0}^\infty\frac{b^kr^k}{4^k(k!)^2},
 \qquad r\in\R.
\]
Its spatial Fourier transform, as in \cite[Section~2]{CL2026} and
\cite[Section~2]{RW2026}, is
\[
\widehat G_{a,m}(t,\xi)
 =e^{-at/2}\frac{\sin(t\sqrt{\xi^2+m^2-a^2/4})}
 {\sqrt{\xi^2+m^2-a^2/4}}.
\]
The expression is interpreted by continuity when the denominator
vanishes, giving $te^{-at/2}$, and by analytic continuation when
the quantity under the square root is negative.
In particular, for every $T>0$,
\begin{equation}\label{eq:greenbound}
 |G_{a,m}(t,x)|
 \le C_T\ind_{\{|x|<t\}},\qquad 0<t\le T.
\end{equation}

A predictable random field $u$ is a mild solution of
\eqref{eq:model} if the stochastic integral below is well defined
and, for every $(t,x)\in \mathbb R_+\times \mathbb R$,
\begin{equation}\label{eq:original-mild}
 u(t,x)=\int_0^t\int_\R
 G_{a,m}(t-s,x-y)F(u(s,y))\,W(ds,dy).
\end{equation}
 
\begin{proposition}\label{prop:regularity}
For globally Lipschitz $F$, equation~\eqref{eq:model} has a unique
predictable mild random-field solution satisfying, for every
$T>0$ and $1\le p<\infty$,
\begin{equation}\label{eq:moments}
 \sup_{0\le t\le T}\sup_{x\in\R}\|u(t,x)\|_p<\infty.
\end{equation}
It also satisfies
\begin{equation}\label{eq:holder}
 \|u(t,x)-u(s,y)\|_p
 \le C_{p,T}\bigl(|t-s|^H+|x-y|^H\bigr),
 \qquad s,t\in[0,T],\quad x,y\in\R.
\end{equation}
\end{proposition}
\begin{proof}
Fix $T>0$ and $p\ge2$, and write $G=G_{a,m}$.
We follow the Picard argument in \cite[Proposition~2.1]{RW2026},
replacing the spatial $L^2$ estimates by Riesz covariance estimates.
For this purpose, recall the interval identity of Liu and Wang
\cite[Lemma~3.2(a)]{LW2025}:
\begin{equation}\label{eq:interval-energy}
 \|\ind_I\|_{\mathcal H}^2=|I|^{2H},
 \qquad I\subset\R\text{ a bounded interval}.
\end{equation}
Together with \eqref{eq:greenbound}, it gives
\begin{equation}\label{eq:green-riesz-bound}
 \bigl\||G(r,\cdot)|\bigr\|_{\mathcal H}^2
 \le C_T\|\ind_{(-r,r)}\|_{\mathcal H}^2
 \le C_T r^{2H},\qquad 0\le r\le T.
\end{equation}

Set $u_0=0$ and define
\begin{equation}\label{eq:riesz-picard}
 u_{n+1}(t,x)=\int_0^t\int_\R
       G(t-s,x-y)F(u_n(s,y))\,W(ds,dy).
\end{equation}
By \eqref{eq:bdg} and \eqref{eq:green-riesz-bound}, each integral
is well defined whenever the preceding iterate has bounded moments.
We use predictable versions of these adapted convolutions.
Their mean-square continuity follows by dominated convergence
from \eqref{eq:fullgreen}--\eqref{eq:greenbound} and the local
integrability of the Riesz kernel; see \cite{Walsh1986}.
For the moment and convergence estimates, put
\[
 M_n(t)=\sup_x\E|u_n(t,x)|^p,\qquad
 D_n(t)=\sup_x\E|u_{n+1}(t,x)-u_n(t,x)|^p.
\]
The linear growth and Lipschitz continuity of $F$, together with
\eqref{eq:bdg}, \eqref{eq:green-riesz-bound}, and H\"older's
inequality in time, yield
\begin{equation}\label{eq:riesz-picard-moments}
\begin{aligned}
 M_{n+1}(t)&\le C_{p,T}\int_0^t(1+M_n(s))\,ds,
       &&n\ge0,\\
 D_n(t)&\le C_{p,T}\int_0^tD_{n-1}(s)\,ds,
       &&n\ge1.
\end{aligned}
\end{equation}
The first inequality gives $M_n(t)\le e^{C_{p,T}t}-1$,
uniformly in $n$. Iterating the second, with
$A=\sup_{t\le T}D_0(t)<\infty$, gives
\[
 D_n(t)\le A\frac{(C_{p,T}t)^n}{n!},\qquad
 \sum_{n=0}^\infty\sup_{t\le T}D_n(t)^{1/p}<\infty.
\]
Hence $u_n$ converges in $L^p(\Omega)$, uniformly in
$(t,x)\in[0,T]\times\R$, to a predictable field $u$ satisfying
\eqref{eq:moments}. By the Lipschitz property of $F$ and
\eqref{eq:bdg}, we may pass to the limit in
\eqref{eq:riesz-picard}, so $u$ is a mild solution.
For two such solutions $u$ and $\widetilde u$, the same estimates
with $p=2$ give
\[
 \sup_x\|u(t,x)-\widetilde u(t,x)\|_2^2
 \le C_T\int_0^t\sup_y
       \|u(s,y)-\widetilde u(s,y)\|_2^2\,ds.
\]
Uniqueness follows from Gronwall's lemma.

It remains to prove \eqref{eq:holder}. Define
\[
 K_{t,x}(r,z)=\ind_{\{0\le r<t\}}G(t-r,x-z),
\]
and write $G(r,z)=\Phi(r,z)\ind_{I_r}(z)$, where
\[
 I_r=(-r,r),\qquad
 \Phi(r,z)=\tfrac12e^{-ar/2}\mathcal B_b(r^2-z^2).
\]
The power series defining $\mathcal B_b$ gives
\[
 \sup_{\substack{0\le r\le T\\ |z|\le T+1}}
 \bigl(|\Phi(r,z)|+|\partial_r\Phi(r,z)|
                    +|\partial_z\Phi(r,z)|\bigr)\le C_T.
\]
For a spatial shift $0<\delta\le1$, put $I_r^\delta=I_r+\delta$.
By the mean value theorem,
\[
 |G(r,z-\delta)-G(r,z)|
 \le C_T\bigl(\delta\ind_{I_r\cap I_r^\delta}(z)
       +\ind_{I_r^\delta\setminus I_r}(z)
       +\ind_{I_r\setminus I_r^\delta}(z)\bigr).
\]
Each difference set is an interval of length at most $\delta$.
Thus \eqref{eq:interval-energy} and the triangle inequality imply
\[
 \bigl\||G(r,\cdot-\delta)-G(r,\cdot)|\bigr\|_{\mathcal H}
 \le C_T\bigl(\delta(2r)^H+2\delta^H\bigr).
\]
Integrating in time and using translation invariance, we obtain
\[
 \int_0^t\bigl\||K_{t,x+\delta}(r,\cdot)
                         -K_{t,x}(r,\cdot)|\bigr\|_{\mathcal H}^2\,dr
 \le C_T\bigl(\delta^2+\delta^{2H}\bigr)
 \le C_T\delta^{2H}.
\]

For temporal shifts, $I_{r+\delta}\setminus I_r$ consists of two
intervals of length $\delta$. The same argument therefore gives
\[
 \bigl\||G(r+\delta,\cdot)-G(r,\cdot)|\bigr\|_{\mathcal H}
 \le C_T\bigl(\delta(2r)^H+2\delta^H\bigr),
 \qquad 0\le r\le r+\delta\le T,
\]
for $0<\delta\le1$. Let $0\le s<t\le T$ and $\delta=t-s\le1$.
Splitting the time integral at $s$ and using
\eqref{eq:green-riesz-bound} on $(s,t)$, we obtain
\[
\begin{aligned}
 &\int_0^t\bigl\||K_{t,x}(r,\cdot)
                    -K_{s,x}(r,\cdot)|\bigr\|_{\mathcal H}^2\,dr\\
 &=\int_0^s
   \bigl\||G(t-r,\cdot)-G(s-r,\cdot)|\bigr\|_{\mathcal H}^2\,dr
   +\int_s^t\bigl\||G(t-r,\cdot)|\bigr\|_{\mathcal H}^2\,dr\\
 &\le C_T\left[\int_0^s
       \bigl(\delta^2(s-r)^{2H}+\delta^{2H}\bigr)\,dr
       +\int_0^\delta r^{2H}\,dr\right]
 \le C_T\delta^{2H}.
\end{aligned}
\]
For increments larger than one, the same bounds follow directly
from \eqref{eq:green-riesz-bound}. Combining the spatial and
temporal estimates gives
\begin{equation}\label{eq:green-riesz-increments}
 \int_0^T\bigl\||K_{t,x}(r,\cdot)
                    -K_{s,y}(r,\cdot)|\bigr\|_{\mathcal H}^2\,dr
 \le C_T\bigl(|t-s|^{2H}+|x-y|^{2H}\bigr).
\end{equation}
Finally, \eqref{eq:original-mild}, \eqref{eq:bdg}, and
\eqref{eq:moments} yield
\[
\begin{aligned}
 \|u(t,x)-u(s,y)\|_p^2
 &\le C_{p,T}\int_0^T
       \bigl\||K_{t,x}(r,\cdot)-K_{s,y}(r,\cdot)|
       \bigr\|_{\mathcal H}^2\,dr\\
 &\le C_{p,T}\bigl(|t-s|^{2H}+|x-y|^{2H}\bigr),
\end{aligned}
\]
which proves \eqref{eq:holder} for $p\ge2$.
The case $1\le p<2$ follows from the second-moment bounds.
Since $p$ can be arbitrarily large, Kolmogorov's theorem gives a
continuous modification, locally H\"older continuous of every
order less than $H$.
\end{proof}

\section{Bi-parameter local linearization}\label{sec:local}

\subsection{The critical kernel and the Klein--Gordon resolvent}
For the mixed-increment analysis, we introduce the critical kernel
\[
 \Gamma_a(t,x)=\frac12e^{-at/2}\ind_{\{t>0,\ |x|<t\}},
\]
which is the fundamental solution of
$\partial_t^2-\partial_x^2+a\partial_t+a^2/4$.
Following Chen and Lee~\cite[Section~3]{CL2026} and
Rang and Wang~\cite[Section~3.1]{RW2026}, we write
\begin{equation}\label{eq:critical-mild}
\begin{split}
 u(t,x)&=u_L(t,x)+u_C(t,x),\\
 u_L(t,x)&=b\int_0^t\int_\R\Gamma_a(t-s,x-y)u(s,y)\,dy\,ds,\\
 u_C(t,x)&=\int_0^t\int_\R\Gamma_a(t-s,x-y)F(u(s,y))\,W(ds,dy).
\end{split}
\end{equation}
The additive solution has the analogous decomposition
$U=U_L+U_C$, with $F\equiv1$.
The stochastic Fubini step in this reduction remains valid for
the present noise. Indeed, \eqref{eq:moments} and
\eqref{eq:green-riesz-bound} imply, for $t\le T$,
\[
\begin{aligned}
 &\int_0^t\int_\R |\Gamma_a(t-s,x-y)|
 \left(\E\int_0^s
 \|G_{a,m}(s-r,y-\cdot)F(u(r,\cdot))\|_{\mathcal H}^2\,dr
 \right)^{1/2}\,dy\,ds\\
  \le &\,  C_T\int_0^t(t-s)s^{H+1/2}\,ds<\infty.
\end{aligned}
\]

\subsection{Characteristic geometry and kernel bounds}
We first consider positive mixed increments $D_h=D_h^+$, with
$0<h\le1$ and $0\le t_z\le T-\sqrt2h$.
For a physical source point $(s,y)$, write
\[
\alpha=\frac{s-y}{\sqrt2},\qquad \zeta=\frac{s+y}{\sqrt2}.
\]
For $z=(\tau,\lambda)$, define the corresponding characteristic cell
\begin{equation}\label{eq:diamond}
 P_z(h)=\{(s,y):\tau<\alpha\le\tau+h,
                     \lambda<\zeta\le\lambda+h\}.
\end{equation}
This characteristic cell is a diamond in physical coordinates,
with vertices
\[
\begin{aligned}
 P_0&=(t_z,x_z),&
 P_1&=(t_z+h/\sqrt2,x_z-h/\sqrt2),\\
 P_2&=(t_z+h/\sqrt2,x_z+h/\sqrt2),&
 P_{12}&=(t_z+\sqrt2h,x_z).
\end{aligned}
\]
In particular, all source times in $P_z(h)$ are larger than $t_z$,
up to an irrelevant boundary. The cell noise
\begin{equation}\label{eq:cell-noise}
 Z_z(h)=\frac12\int\ind_{P_z(h)}(s,y)\,W(ds,dy)
\end{equation}
is independent of $\mathcal F_{t_z}$, by temporal independence of
the Gaussian noise.

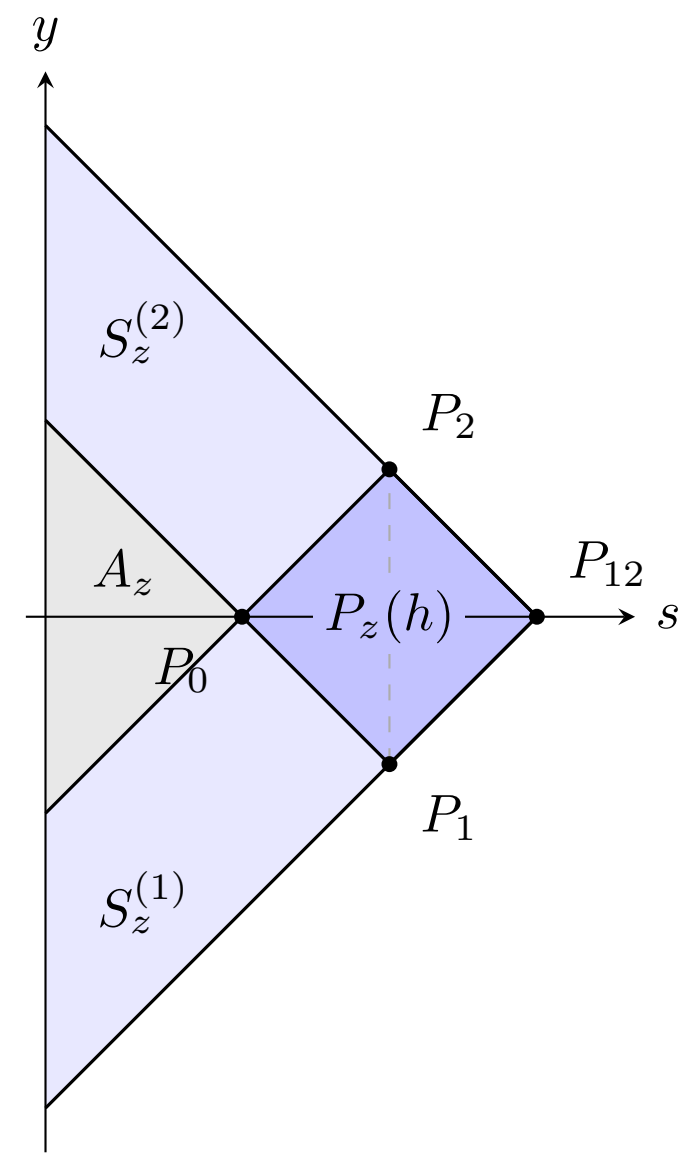
\begin{figure}[tb]
\centering
\begin{tikzpicture}[x=0.65cm,y=0.65cm,>=stealth,
                    line join=round,font=\small]
 \fill[gray!18] (0,-2)--(0,2)--(2,0)--cycle;
 \fill[blue!9] (0,-5)--(0,-2)--(2,0)--(3.5,-1.5)--cycle;
 \fill[blue!9] (0,2)--(0,5)--(3.5,1.5)--(2,0)--cycle;
 \fill[blue!24] (2,0)--(3.5,-1.5)--(5,0)--(3.5,1.5)--cycle;
 \draw[thin,gray!70] (0,-5)--(3.5,-1.5)--(0,2);
 \draw[thin,gray!70] (0,-2)--(3.5,1.5)--(0,5);
 \draw[semithick] (0,-5)--(5,0)--(0,5);
 \draw[semithick] (0,-2)--(2,0)--(0,2);
 \draw[semithick] (2,0)--(3.5,-1.5)--(5,0)--(3.5,1.5)--cycle;
 \draw[dashed,gray!65] (3.5,-1.5)--(3.5,1.5);
 \draw[->] (0,-5.45)--(0,5.55) node[above]{$y$};
 \draw[->] (-.2,0)--(6,0) node[right]{$s$};
 \node at (.8,.45){$A_z$};
 \node at (1,-2.9){$S_z^{(1)}$};
 \node at (1,2.9){$S_z^{(2)}$};
 \node[fill=blue!24,inner sep=2pt] at (3.5,0){$P_z(h)$};
 \fill (2,0) circle (1.5pt) node[below left=2pt]{$P_0$};
 \fill (3.5,-1.5) circle (1.5pt) node[below right=2pt]{$P_1$};
 \fill (3.5,1.5) circle (1.5pt) node[above right=2pt]{$P_2$};
 \fill (5,0) circle (1.5pt) node[above right=2pt]{$P_{12}$};
\end{tikzpicture}
\caption{The old cone $A_z$, characteristic strips $S_z^{(1)}$ and
$S_z^{(2)}$, and new diamond $P_z(h)$ in physical coordinates.}
\label{fig:geometry}
\end{figure}

\begin{lemma}[{\cite[Lemma~3.1, equation~(3.1)]{LW2025}}]\label{lem:cell-variance}
For every pair $(z,h)$ admissible for the positive orientation,
\begin{equation}\label{eq:cell-variance}
 \E Z_z(h)^2=\kappa_H h^{2H+1},
\end{equation}
where $\kappa_H$ is defined in \eqref{eq:kappa-intro}.
Consequently, $\|Z_z(h)\|_p=C_{p,H}h^{H+1/2}$ for every
$1\le p<\infty$.
\end{lemma}
Define the critical mixed kernel by
\begin{equation}\label{eq:mixed-kernel}
\begin{aligned}
 K_{z,h}(s,y)={}&\Gamma_a(t_z+\sqrt2h-s,x_z-y)-\Gamma_a(t_z+h/\sqrt2-s,x_z-h/\sqrt2-y)\\
 &-\Gamma_a(t_z+h/\sqrt2-s,x_z+h/\sqrt2-y)
   +\Gamma_a(t_z-s,x_z-y).
\end{aligned}
\end{equation}
The remaining regions needed for its decomposition are
\[
\begin{aligned}
 A_z&=\{(s,y):s\ge0,\ \alpha\le\tau,\ \zeta\le\lambda\},\\
 S_z^{(1)}&=\{(s,y):s\ge0,\ \tau<\alpha\le\tau+h,\ \zeta\le\lambda\},\\
 S_z^{(2)}&=\{(s,y):s\ge0,\ \alpha\le\tau,\ \lambda<\zeta\le\lambda+h\}.
\end{aligned}
\]
Figure~\ref{fig:geometry} illustrates these regions.

\begin{lemma}\label{lem:kernel-bounds}
For $0<h\le1$ and $0\le t_z\le t_z+\sqrt2h\le T$,
\begin{align}
 \int_0^\infty\int_\R |K_{z,h}(s,y)|\,dy\,ds
    &\le C_T h^2,\label{eq:mixed-lone}\\
 \int_0^\infty
 \bigl\||K_{z,h}(s,\cdot)-\tfrac12\ind_{P_z(h)}(s,\cdot)|\bigr\|_{\mathcal H}^2\,ds
    &\le C_T h^{2H+2},\label{eq:tail-energy}\\
 \int_0^\infty\bigl\||K_{z,h}(s,\cdot)|\bigr\|_{\mathcal H}^2\,ds
    &\le C_T h^{2H+1}.\label{eq:mixed-energy}
\end{align}
\end{lemma}
\begin{proof}
Estimate \eqref{eq:mixed-lone} follows from
\cite[Lemma~3.1]{RW2026} with $p=1$, $\varepsilon=h$, and
$(t,x)=(t_z+\sqrt2h,x_z)$. We use the same decomposition to
prove the two Riesz norm estimates. Set
\[
 q_h=e^{-ah/(2\sqrt2)},\qquad
 w_z(s)=\tfrac12e^{-a(t_z-s)/2}.
\]
In the present notation, the decomposition reads
\begin{equation}\label{eq:kernel-decomposition}
 K_{z,h}=w_z\left[(q_h-1)^2\ind_{A_z}
    +q_h(q_h-1)\left(\ind_{S_z^{(1)}}+\ind_{S_z^{(2)}}\right)
    +q_h^2\ind_{P_z(h)}\right].
\end{equation}
On $0\le s\le T$, we have $|w_z(s)|+|q_h|\le C_T$ and
$|q_h-1|\le C_Th$. On the diamond $P_z(h)$,
$0\le t_z+\sqrt2h-s\le\sqrt2h$, so
\[
 w_z(s)q_h^2=\tfrac12e^{-a(t_z+\sqrt2h-s)/2}
           =\tfrac12+O_T(h).
\]
Consequently, almost everywhere,
\[
 \bigl|K_{z,h}-\tfrac12\ind_{P_z(h)}\bigr|
 \le C_T\left[h^2\ind_{A_z}
       +h\left(\ind_{S_z^{(1)}}+\ind_{S_z^{(2)}}\right)
       +h\ind_{P_z(h)}\right].
\]

The spatial sections of $A_z$ have length $2(t_z-s)$ for
$0\le s\le t_z$, while those of each strip are intervals of
length at most $\sqrt2h$. Since all regions lie in
$[0,T]\times\R$, \eqref{eq:interval-energy} and
Lemma~\ref{lem:cell-variance} give
\[
\begin{aligned}
 \int_0^\infty\|\ind_{A_z}(s,\cdot)\|_{\mathcal H}^2\,ds
 &=\int_0^{t_z}[2(t_z-s)]^{2H}\,ds\le C_T,\\
 \int_0^\infty\|\ind_{S_z^{(j)}}(s,\cdot)\|_{\mathcal H}^2\,ds
 &\le T(\sqrt2h)^{2H},\qquad j=1,2,\\
 \int_0^\infty\|\ind_{P_z(h)}(s,\cdot)\|_{\mathcal H}^2\,ds
 &=4\kappa_Hh^{2H+1}.
\end{aligned}
\]
Combining the pointwise bound with these estimates, using
positivity of the Riesz covariance kernel and the triangle
inequality in $L^2(\R_+;\mathcal H)$, yields
\[
\begin{aligned}
 &\int_0^\infty
 \bigl\||K_{z,h}(s,\cdot)-\tfrac12\ind_{P_z(h)}(s,\cdot)|
       \bigr\|_{\mathcal H}^2\,ds\\
  \le &\, C_T\bigl(h^4+h^{2H+2}+h^{2H+3}\bigr)
 \le C_Th^{2H+2},
\end{aligned}
\]
since $H<1$ and $h\le1$. This proves \eqref{eq:tail-energy}.
Adding back $\tfrac12\ind_{P_z(h)}$ and applying the same
inequalities gives
\[
 \int_0^\infty\bigl\||K_{z,h}(s,\cdot)|\bigr\|_{\mathcal H}^2\,ds
 \le C_T\bigl(h^{2H+2}+h^{2H+1}\bigr)
 \le C_Th^{2H+1},
\]
which proves \eqref{eq:mixed-energy}.
\end{proof}

\subsection{Local approximation}
For the drift components in \eqref{eq:critical-mild}, write
\[
 v_L(z)=u_L(t_z,x_z),\qquad V_L(z)=U_L(t_z,x_z).
\]

The kernel bounds give the following approximation by the cell
noise $Z_z(h)$.

\begin{proposition}\label{prop:cell-approx}
Under the assumptions of Theorem~\ref{thm:linearization}, for every
$1\le p<\infty$ and $(z,h)$ admissible for the positive orientation
with $0<h\le1$,
\begin{align}
 \|D_hv_L(z)\|_p+\|D_hV_L(z)\|_p&\le C_{p,T}h^2,
       \label{eq:drift-mixed}\\
 \|D_hv(z)-F(v(z))Z_z(h)\|_p&\le C_{p,T}h^{H+1},
       \label{eq:nonlinear-cell}\\
 \|D_hV(z)-Z_z(h)\|_p&\le C_{p,T}h^{H+1}.
       \label{eq:linear-cell}
\end{align}
Moreover,
\begin{equation}\label{eq:increment-size}
 \|D_hv(z)\|_p+\|D_hV(z)\|_p\le C_{p,T}h^{H+1/2}.
\end{equation}
\end{proposition}
\begin{proof}
It suffices to consider $p\ge2$. The drift estimate
\eqref{eq:drift-mixed} is the deterministic convolution estimate
of \cite[Proposition~3.1]{RW2026}. Its proof uses only the
$L^1$ kernel bound \eqref{eq:mixed-lone} and uniform moment bounds,
so it applies here by Proposition~\ref{prop:regularity}.

Put $L_{z,h}=K_{z,h}-\tfrac12\ind_{P_z(h)}$. We have
\begin{equation}\label{eq:freeze-decomposition}
\begin{aligned}
 D_hv(z)-F(v(z))Z_z(h)
 ={}&D_hv_L(z)
   +\int L_{z,h}(s,y)F(u(s,y))\,W(ds,dy)\\
 &+\frac12\int_{P_z(h)}[F(u(s,y))-F(u(t_z,x_z))]\,W(ds,dy).
\end{aligned}
\end{equation}
The frozen coefficient is $\mathcal F_{t_z}$-measurable and the
cell lies after $t_z$, so the last integrand is predictable.
By \eqref{eq:bdg}, \eqref{eq:moments}, and \eqref{eq:tail-energy},
the second term has $L^p$ norm at most $C_{p,T}h^{H+1}$.
The cell-freezing estimate in the proof of
\cite[Theorem~1.1]{LW2025} uses only the Lipschitz property of $F$
and the moment and increment bounds
\eqref{eq:moments}--\eqref{eq:holder}. The same argument bounds
the $L^p$ norm of the last term by $C_{p,T}h^{2H+1/2}$.
Consequently,
\begin{equation}\label{eq:three-errors}
 \|D_hv(z)-F(v(z))Z_z(h)\|_p
 \le C_{p,T}\bigl(h^2+h^{H+1}+h^{2H+1/2}\bigr)
 \le C_{p,T}h^{H+1}.
\end{equation}
This proves \eqref{eq:nonlinear-cell}, and $F\equiv1$ gives
\eqref{eq:linear-cell}. The increment bound
\eqref{eq:increment-size} follows from these estimates,
Lemma~\ref{lem:cell-variance}, and \eqref{eq:moments}.
\end{proof}

\begin{proof}[Proof of Theorem~\ref{thm:linearization}]
For the positive orientation, add and subtract $F(v(z))Z_z(h)$.
Proposition~\ref{prop:cell-approx} and H\"older's inequality give
\[
\begin{aligned}
 \|R_h^+(z)\|_p
  \le&\, \|D_hv(z)-F(v(z))Z_z(h)\|_p+\|F(v(z))\|_{2p}\|D_hV(z)-Z_z(h)\|_{2p}\\
  \le &\,  C_{p,T}h^{H+1}.
\end{aligned}
\]
For the negative orientation, suppose
$h/\sqrt2\le t_z\le T-h/\sqrt2$ and set $z_-=z-(h,0)$.
Then $(z_-,h)$ is admissible for the positive orientation and
$D_h^-f(z)=-D_hf(z_-)$. Consequently,
\[
 R_h^-(z)=-R_h^+(z_-)
       +[F(v(z))-F(v(z_-))]D_hV(z_-).
\]
By \eqref{eq:holder} and \eqref{eq:increment-size},
\[
 \|[F(v(z))-F(v(z_-))]D_hV(z_-)\|_p
 \le C_{p,T} h^{2H+1/2}.
\] 
Combining this with the estimate for the positive orientation
proves \eqref{eq:main-linearization} for both signs.
\end{proof}

\section{Quadratic variations and parameter estimation}
\label{qv:section}

\subsection{Gaussian cell increments and weighted squares}
 
Set
\[
 Z_{ij}=Z_{z_{ij}}\left(N^{-1}\right),\qquad
 \sigma_N^2=\E\left[Z_{ij}^2\right]=\kappa_H N^{-2H-1}.
\]
The following $L^p$ estimate yields the almost-sure convergence
needed for strong consistency.

\begin{lemma}[Predictable weighted squares]\label{qv:weighted}
Fix $p\ge2$. Suppose that $a_{ij}$ is
$\mathcal F_{(i+j)/(\sqrt2N)}$-measurable and
$\sup_{N,i,j}\norm{a_{ij}}_{L^p(\Omega)}\le K_p<\infty$.
Then
\begin{equation}\label{qv:weighted-bound}
 \left\|N^{2H-1}\sum_{i,j=0}^{N-1}
       a_{ij}\left(Z_{ij}^2-\sigma_N^2\right)\right\|_{L^p(\Omega)}
 \le C_{p, H} K_p N^{-1/2}.
\end{equation} 
\end{lemma}

\begin{proof}
We use the physical-time grouping in the proof of
\cite[Proposition~1.3]{LW2025}, with the martingale inequality
in $L^p(\Omega)$. Put $t_k=k/(\sqrt2N)$ and
\[
 D_k=\sum_{\substack{0\le i,j<N\\i+j=k}}
        a_{ij}(Z_{ij}^2-\sigma_N^2),
 \qquad
 n_k=\#\{(i,j):0\le i,j<N,\ i+j=k\}.
\]
The cells with $i+j=k$ lie between physical times $t_k$ and
$t_{k+2}$. Since the noise is white in time, the vector
$\{Z_{ij}:i+j=k\}$ is independent of $\mathcal F_{t_k}$ and
is $\mathcal F_{t_{k+2}}$-measurable. Thus $D_k$ is
$\mathcal F_{t_{k+2}}$-measurable and
$\E[D_k\mid\mathcal F_{t_k}]=0$.
For each $e\in\{0,1\}$, the variables $D_{2r+e}$ are therefore
martingale differences between
$\mathcal F_{t_{2r+e}}$ and $\mathcal F_{t_{2r+e+2}}$.
Independence of each weight from its cell noise and Gaussian
moment bounds give
$\|D_k\|_p\le C_pK_p n_k\sigma_N^2$.
Applying Burkholder's square-function inequality to the even
and odd martingale sums, followed by Minkowski's inequality in
$L^{p/2}(\Omega)$, gives, since $p\ge2$,
\[
\begin{aligned}
 \left\|\sum_{k=0}^{2N-2}D_k\right\|_p
 &\le C_p\sum_{e=0}^1
    \left\|\left(
       \sum_{\substack{0\le k\le2N-2\\ k\equiv e\pmod{2}}}D_k^2
    \right)^{1/2}\right\|_p\\
 &\le C_p\left(\sum_{k=0}^{2N-2}\|D_k\|_p^2\right)^{1/2}
 \le C_pK_p\sigma_N^2N^{3/2},
\end{aligned}
\]
since $\sum_k n_k^2\le N^3$.
Multiplying by $N^{2H-1}$ and using
$N^{2H-1}\sigma_N^2=\kappa_HN^{-2}$ proves the claim.
\end{proof}

\subsection{Convergence of the nonlinear quadratic variation}
We combine Proposition~\ref{prop:cell-approx} and
Lemma~\ref{qv:weighted} in the decomposition used in
\cite[Proposition~1.3]{LW2025}.

\begin{proof}[Proof of Theorem~\ref{thm:variation}]
It suffices to consider $p\ge2$. Write $f_{ij}=F(v(z_{ij}))$ and set
\[
A=\int_0^1\int_0^1 F^2(v(\tau,\lambda))\,d\tau\,d\lambda.
\]
Decompose
\begin{align*}
 Q_N(v)-\kappa_H A
 ={}&N^{2H-1}\sum_{i,j}
        \bigl((\Delta_{ij}^N v)^2-f_{ij}^2Z_{ij}^2\bigr)\\
 &+N^{2H-1}\sum_{i,j}f_{ij}^2\left(Z_{ij}^2-\sigma_N^2\right)
   +\kappa_H\left(N^{-2}\sum_{i,j}f_{ij}^2-A\right)\\
 =:{}&A_N+B_N+C_N.
\end{align*}

Proposition~\ref{prop:cell-approx}, Gaussian moment bounds, and
H\"older's inequality give
\begin{align*}
 \norm{A_N}_{L^p(\Omega)}
 &\le N^{2H-1}\sum_{i,j}
     \norm{\Delta_{ij}^N v-f_{ij}Z_{ij}}_{L^{2p}(\Omega)}
     \norm{\Delta_{ij}^N v+f_{ij}Z_{ij}}_{L^{2p}(\Omega)}\\
 &\le C_p N^{2H-1}N^2
            N^{-H-1}N^{-H-1/2}
 =C_p N^{-1/2}.
\end{align*}

The weights $f_{ij}^2$ are measurable with respect to
$\mathcal F_{(i+j)/(\sqrt2N)}$ and satisfy
$\sup_{N,i,j}\|f_{ij}^2\|_p<\infty$ by \eqref{eq:moments}.
Lemma~\ref{qv:weighted}, with $a_{ij}=f_{ij}^2$, therefore gives
\[
\norm{B_N}_{L^p(\Omega)}\le C_p N^{-1/2}.
\]

For $C_N$, the Riemann-sum estimate in the proof of
\cite[Proposition~1.3]{LW2025} applies with H\"older's inequality
in $L^{2p}$ and \eqref{eq:moments}--\eqref{eq:holder}, giving
$\|C_N\|_p\le C_pN^{-H}$.
Because $H>1/2$, these estimates prove \eqref{eq:main-variation}.
For any fixed $p>2$ and $\eta>0$, Markov's inequality gives
\[
\sum_{N=1}^\infty
 \mathbb P\bigl(|Q_N(v)-\kappa_H A|>\eta\bigr)
 \le C_p\eta^{-p}\sum_{N=1}^\infty N^{-p/2}<\infty.
\]
The Borel--Cantelli lemma, applied to the countable sequence
$\eta=1,1/2,1/3,\ldots$, proves the almost-sure assertion.
\end{proof}

\subsection{Diffusion-parameter estimation}\label{sec:estimation}
We apply Theorem~\ref{thm:variation} to the parametric equation
\begin{equation}\label{qv:parameter-model}
 \begin{cases}
  (\partial_t^2-\partial_x^2+a\partial_t+m^2)u_\theta(t,x)
     =\theta F(u_\theta(t,x))\dot W(t,x),&t>0,\ x\in\R,\\
  u_\theta(0,x)=\partial_tu_\theta(0,x)=0,&x\in\R.
 \end{cases}
\end{equation}
We first assume that $a,m,H$, and $F$ are known and estimate
$\theta>0$; Section~\ref{sec:hurst-estimation} treats unknown $H$.
As before, $v_\theta(z)=u_\theta(t_z,x_z)$. We impose
\begin{equation}\label{qv:nondegeneracy}
 F(0)\ne0.
\end{equation}
For zero initial data, this condition excludes the uninformative
case: if $F(0)=0$, uniqueness gives $u_\theta\equiv0$ for every
$\theta$.

Continuity and \eqref{qv:nondegeneracy} imply
\[
 A_\theta:=\int_{[0,1]^2}F^2(v_\theta(z))\,dz>0
 \quad\text{almost surely}.
\]
Theorem~\ref{thm:variation}, applied with coefficient $\theta F$,
gives $Q_N(v_\theta)\to\theta^2\kappa_H A_\theta$ almost surely.
We assume that $u_\theta$ is observed on the characteristic grid
\begin{equation}\label{eq:observation-grid}
 \left\{u_\theta\left(\frac{i+j}{\sqrt2N},
                       \frac{j-i}{\sqrt2N}\right):
                 0\le i,j\le N\right\}.
\end{equation}
Replacing the integral in the limit by its Riemann sum suggests
the estimator
\begin{align}
 \widehat\theta_N
 =\left(
   \frac{N^{2H+1}\displaystyle\sum_{i,j=0}^{N-1}
                          \left(\Delta_{ij}^N v_\theta\right)^2}
   {\kappa_H\displaystyle\sum_{i,j=0}^{N-1}
                          F^2(v_\theta(z_{ij}))}
   \right)^{1/2}.\label{eq:estimator-sums}
\end{align}
The denominator is positive for every $N$, since the sum includes
$F^2(v_\theta(0,0))=F^2(0)>0$.
The restriction $\theta>0$ fixes the sign, since quadratic variations
identify $\theta^2$.

\begin{corollary}\label{qv:consistency}
Under \eqref{qv:nondegeneracy}, the estimator
\eqref{eq:estimator-sums} satisfies
\begin{equation}\label{eq:theta-known-rate}
 \widehat\theta_N\longrightarrow\theta
 \quad\text{almost surely},\qquad
 \widehat\theta_N-\theta=O_{\mathbb P}(N^{-1/2}).
\end{equation}
\end{corollary}
\begin{proof}
Put
\[
 A_{\theta,N}=N^{-2}\sum_{i,j=0}^{N-1}F^2(v_\theta(z_{ij})).
\]
By continuity, $A_{\theta,N}\to A_\theta>0$ almost surely.
The ratio argument in \cite[Corollary~4.1]{LW2025} and
\cite[Corollary~4.1]{RW2026}, now using
Theorem~\ref{thm:variation}, gives the almost-sure convergence.
For the rate, the Riemann-sum estimate in the proof of
Theorem~\ref{thm:variation} gives
\[
 \|A_{\theta,N}-A_\theta\|_p\le C_pN^{-H},
 \qquad 1\le p<\infty.
\]
Since $H>1/2$, Theorem~\ref{thm:variation} implies
\[
 Q_N(v_\theta)-\theta^2\kappa_H A_{\theta,N}
 =O_{\mathbb P}(N^{-1/2}).
\]
Moreover, $A_{\theta,N}^{-1}\to A_\theta^{-1}<\infty$
almost surely, so $A_{\theta,N}^{-1}=O_{\mathbb P}(1)$.
Therefore,
\[
 \widehat\theta_N^2-\theta^2
 =\frac{Q_N(v_\theta)-\theta^2\kappa_H A_{\theta,N}}
        {\kappa_H A_{\theta,N}}
 =O_{\mathbb P}(N^{-1/2}).
\]
The claimed rate follows from
$|\widehat\theta_N-\theta|
 \le\theta^{-1}|\widehat\theta_N^2-\theta^2|$.
\end{proof}

\subsection{Estimation of the Hurst parameter}\label{sec:hurst-estimation}
We now treat both $H\in(1/2,1)$ and $\theta>0$ as unknown
parameters in \eqref{qv:parameter-model}, retaining assumption
\eqref{qv:nondegeneracy}. Define the quadratic variation
\begin{equation}\label{eq:unnormalized-variation}
 S_N(v_\theta)=\sum_{i,j=0}^{N-1}\left(\Delta_{ij}^N v_\theta\right)^2,
 \qquad Q_N(v_\theta)=N^{2H-1}S_N(v_\theta).
\end{equation}
Theorem~\ref{thm:variation} implies that $S_N$ is asymptotically
$N^{1-2H}$ times the same positive random factor at every resolution.
This suggests
\begin{equation}\label{eq:hurst-estimator}
 \widehat H_N
 =\frac12+\frac{1}{2\log2}
       \log\!\left(\frac{S_N(v_\theta)}{S_{2N}(v_\theta)}\right)
 \quad\text{if }S_N(v_\theta)S_{2N}(v_\theta)>0.
\end{equation}
We set $\widehat H_N=3/4$ otherwise. Both sums use the same
sample path on $[0,1]^2$ and are available from the observation
grid \eqref{eq:observation-grid} with $N$ replaced by $2N$.

\begin{proposition}\label{prop:hurst-estimation}
For each fixed $H\in(1/2,1)$ and $\theta>0$, under the assumptions
of Theorem~\ref{thm:variation} and \eqref{qv:nondegeneracy},
\begin{equation}\label{eq:hurst-consistency}
 \widehat H_N\longrightarrow H\quad\text{almost surely},
 \qquad
 \widehat H_N-H=O_{\mathbb P}(N^{-1/2}).
\end{equation}
Moreover, for every $\varepsilon\in(0,1/2)$,
\begin{equation}\label{eq:hurst-as-rate}
 |\widehat H_N-H|=O\left(N^{-1/2+\varepsilon}\right)
 \quad\text{almost surely}.
\end{equation}
\end{proposition}
\begin{proof}
Put
\[
 L=\theta^2\kappa_H A_\theta,
 \qquad E_N=Q_N(v_\theta)-L.
\]
By the preceding subsection, $L>0$ almost surely.
Theorem~\ref{thm:variation} yields
\begin{equation}\label{eq:hurst-qv-error}
\begin{gathered}
 \|E_N\|_p\le C_pN^{-1/2},\qquad 1\le p<\infty,\\
 Q_N(v_\theta)\longrightarrow L\quad\text{almost surely}.
\end{gathered}
\end{equation}
Thus both sums in \eqref{eq:hurst-estimator} are positive for all
sufficiently large $N$, almost surely. On this event,
\[
 \widehat H_N-H
 =\frac{\log Q_N(v_\theta)-\log Q_{2N}(v_\theta)}{2\log2}
 \longrightarrow0.
\]
Whenever $|E_N|\vee|E_{2N}|\le L/2$, the mean value theorem gives
\begin{equation}\label{eq:hurst-log-bound}
 |\widehat H_N-H|
 \le\frac{|E_N|+|E_{2N}|}{L\log2}.
\end{equation}
Localizing on $\{L\ge\eta\}$, using \eqref{eq:hurst-qv-error},
and then letting $\eta\downarrow0$ proves the bound in probability.
For \eqref{eq:hurst-as-rate}, choose $p>1/\varepsilon$.
Markov's inequality gives
\[
 \sum_{N=1}^\infty
 \mathbb P\bigl(|E_N|>N^{-1/2+\varepsilon}\bigr)
 \le C_p\sum_{N=1}^\infty N^{-p\varepsilon}<\infty.
\]
The Borel--Cantelli lemma and \eqref{eq:hurst-log-bound} complete
the proof.
\end{proof}

The almost-sure rate also permits estimation of $\theta$ when $H$
is unknown. Define
\[
 \overline H_N:=\min\{1,\max\{1/2,\widehat H_N\}\}
\]
and set
\begin{equation}\label{eq:joint-diffusion-estimator}
 \widetilde\theta_N
:=\left(
 \frac{N^{2\overline H_N+1}S_N(v_\theta)}
 {\kappa_{\overline H_N}\displaystyle\sum_{i,j=0}^{N-1}
      F^2(v_\theta(z_{ij}))}
 \right)^{1/2}.
\end{equation}
Here $\kappa_H$ is extended continuously to $[1/2,1]$ by
\eqref{eq:kappa-intro}; it is positive throughout this interval.

\begin{corollary}[Joint consistency and rates]\label{cor:joint-consistency}
Under the assumptions of Proposition~\ref{prop:hurst-estimation},
with $F$ known,
\[
 (\widehat H_N,\widetilde\theta_N)
 \longrightarrow(H,\theta)\quad\text{almost surely}.
\]
Moreover,
\begin{equation}\label{eq:joint-rates}
 \widehat H_N-H=O_{\mathbb P}(N^{-1/2}),\qquad
 \widetilde\theta_N-\theta
 =O_{\mathbb P}(N^{-1/2}\log N).
\end{equation}
\end{corollary}
\begin{proof}
Let $\widehat\theta_N$ be the estimator \eqref{eq:estimator-sums} evaluated
at the true $H$. Comparing \eqref{eq:joint-diffusion-estimator}
with \eqref{eq:estimator-sums}, we obtain
\[
 \widetilde\theta_N=\widehat\theta_N e^{\ell_N},\qquad
 \ell_N=(\overline H_N-H)\log N
       +\tfrac12\bigl(\log\kappa_H-\log\kappa_{\overline H_N}\bigr).
\]
The function $r\mapsto\log\kappa_r$ is continuously
differentiable on $[1/2,1]$, and the projection satisfies
$|\overline H_N-H|\le|\widehat H_N-H|$. Consequently,
\begin{equation}\label{eq:joint-log-error}
 |\ell_N|\le C(1+\log N)|\widehat H_N-H|.
\end{equation}
By Corollary~\ref{qv:consistency} and
Proposition~\ref{prop:hurst-estimation}, there is an event
$\Omega_0$ of probability one on which
\[
 \widehat\theta_N\to\theta,\qquad
 \widehat H_N\to H,\qquad
 |\widehat H_N-H|=O(N^{-1/4}).
\]
On this common event, \eqref{eq:joint-log-error} gives
$\ell_N\to0$, so $\widetilde\theta_N=\widehat\theta_Ne^{\ell_N}
\to\theta$. Thus
\[
 |\widehat H_N-H|+|\widetilde\theta_N-\theta|\longrightarrow0
 \quad\text{on }\Omega_0,
\]
which proves joint almost-sure convergence.

For the rates, Proposition~\ref{prop:hurst-estimation} and
\eqref{eq:joint-log-error} imply
\[
 \ell_N=O_{\mathbb P}(N^{-1/2}\log N)=o_{\mathbb P}(1).
\]
Taylor's formula gives
$e^{\ell_N}-1=O_{\mathbb P}(N^{-1/2}\log N)$.
Using \eqref{eq:theta-known-rate} and
$\widehat\theta_N=O_{\mathbb P}(1)$ in
\[
 \widetilde\theta_N-\theta
 =\widehat\theta_N-\theta
   +\widehat\theta_N(e^{\ell_N}-1)
\]
proves the second bound in \eqref{eq:joint-rates}.
\end{proof}
The logarithmic factor in \eqref{eq:joint-rates} comes from
estimating the exponent of $N$ in the diffusion estimator.
These bounds do not require inverse-moment assumptions on
$A_\theta$, since $A_{\theta,N}^{-1}=O_{\mathbb P}(1)$ suffices.

\section{Numerical illustration}\label{sec:simulation}
We illustrate joint estimation of $(\theta,H)$ from
characteristic-grid observations.

\subsection{Simulation setup}
We consider
\[
 \theta=1,\qquad a=m=1,\qquad F(r)=1+\tfrac12\tanh r,
 \qquad H\in\{0.875,0.750,0.625\}.
\]
This coefficient is globally Lipschitz and satisfies $F(0)=1$,
so \eqref{qv:nondegeneracy} holds. For each $H$, we generate
$B=50$ independent noise realizations at simulation resolution
$n=512$, and take $N=8,16,32,64,128$ in the estimators.

To discretize the equation, put $w(t,x)=e^{at/2}u_\theta(t,x)$
and $b=a^2/4-m^2=-3/4$. Set
$\delta=(\sqrt2n)^{-1}$, $t_k=k\delta$, and $x_l=l\delta$.
Let $w_{k,l}$ approximate $w(t_k,x_l)$ and recover
$u_{k,l}=e^{-at_k/2}w_{k,l}$. We compute these values for integer
indices satisfying
\[
 0\le k\le2n,\qquad k+l\ \text{even},\qquad k+|l|\le2n.
\]
This domain contains the backward cones of all observation points,
so no artificial boundary values are required.

The noise integrals are defined over the triangles
\[
\begin{aligned}
 T_{k,l}^{\uparrow}
 &=\{(s,y):k\delta<s<(k+1)\delta,\ |y-l\delta|<s-k\delta\},\\
 T_{k,l}^{\downarrow}
 &=\{(s,y):k\delta<s<(k+1)\delta,\ |y-l\delta|<(k+1)\delta-s\}.
\end{aligned}
\]
Set $\xi_{k,l}=W(T_{k,l}^{\uparrow})$ for even $k+l$ and
$\xi_{k,l}=W(T_{k,l}^{\downarrow})$ otherwise. Noise integrals
over different time strips are independent. Within each strip,
the Gaussian vector
has the Toeplitz covariance obtained by integrating
\eqref{eq:noise-intro} over pairs of triangles. We sample the
noise vectors by circulant embedding and the fast Fourier transform (FFT); see
\cite{DaviesBryant2013}.

On the prescribed domain, initialize $w_{0,l}=0$ for even $l$ and
$w_{1,l}=\theta F(0)\xi_{0,l}/2$ for odd $l$.
For $0\le k\le2n-2$, $k+l$ even, and $k+2+|l|\le2n$, use
\begin{equation}\label{sim:scheme}
\begin{aligned}
 w_{k+2,l}={}&w_{k+1,l-1}+w_{k+1,l+1}-w_{k,l}
                +b\delta^2w_{k,l}\\
 &+\frac{\theta}{2}e^{at_k/2}F(u_{k,l})
          (\xi_{k,l}+\xi_{k+1,l}).
\end{aligned}
\end{equation}
The drift term and the diffusion coefficient are evaluated at $(t_k,x_l)$,
the earliest vertex of the characteristic cell.

We choose $n$ to be an integer multiple of $2N$. The observations
on the finer characteristic grid are read from the simulation as
\[
 v_\theta\left(\frac{i}{2N},\frac{j}{2N}\right)
 \approx u_{k,l},\qquad
 k=\frac{(i+j)n}{2N},\quad l=\frac{(j-i)n}{2N},
 \qquad 0\le i,j\le2N.
\]
The coarser grid consists of the points with both $i$ and $j$ even.
We compute $S_N$ and $S_{2N}$ from these nested observations and
apply \eqref{eq:hurst-estimator} and
\eqref{eq:joint-diffusion-estimator}. For each estimator, we report
the empirical mean and root mean squared error (RMSE).

\subsection{Numerical results}
Table~\ref{tab:simulation-joint} summarizes the estimates
$\widehat H_N$ and $\widetilde\theta_N$; the projection
$\overline H_N$ is used only in $\widetilde\theta_N$.

\begin{table}[tb]
\centering
\caption{Joint estimation of $(\theta,H)$ with $\theta=a=m=1$,
$B=50$, and $n=512$.}
\label{tab:simulation-joint}
\small
\setlength{\tabcolsep}{4pt}
\begin{tabular}{rrrrrr}
\hline
& & \multicolumn{2}{c}{$\widehat H_N$}
  & \multicolumn{2}{c}{$\widetilde\theta_N$}\\
$H$ & $N$ & Mean & RMSE & Mean & RMSE\\
\hline
0.875 & 8 & 0.8265 & 0.1220 & 0.8770 & 0.3093\\
0.875 & 16 & 0.8402 & 0.0791 & 0.8980 & 0.2487\\
0.875 & 32 & 0.8606 & 0.0405 & 0.9485 & 0.1818\\
0.875 & 64 & 0.8697 & 0.0251 & 0.9790 & 0.1278\\
0.875 & 128 & 0.8755 & 0.0164 & 1.0036 & 0.0996\\
\hline
0.750 & 8 & 0.7158 & 0.1203 & 0.9493 & 0.3334\\
0.750 & 16 & 0.7409 & 0.0534 & 0.9788 & 0.1833\\
0.750 & 32 & 0.7485 & 0.0232 & 0.9939 & 0.0998\\
0.750 & 64 & 0.7471 & 0.0162 & 0.9864 & 0.0779\\
0.750 & 128 & 0.7492 & 0.0069 & 0.9947 & 0.0383\\
\hline
0.625 & 8 & 0.6141 & 0.1017 & 0.9997 & 0.2670\\
0.625 & 16 & 0.6124 & 0.0558 & 0.9619 & 0.1871\\
0.625 & 32 & 0.6160 & 0.0303 & 0.9614 & 0.1298\\
0.625 & 64 & 0.6183 & 0.0155 & 0.9657 & 0.0766\\
0.625 & 128 & 0.6213 & 0.0083 & 0.9785 & 0.0464\\
\hline
\end{tabular}
\end{table}

\begin{figure}[tb]
\centering
\begin{tikzpicture}
\begin{groupplot}[group style={group size=2 by 1,horizontal sep=1.8cm},
 width=0.43\textwidth,height=5.2cm,
 xmode=log,log basis x=2,ymode=log,
 xtick={8,16,32,64,128},xticklabels={8,16,32,64,128},
 xlabel={$N$},ylabel={RMSE},
 tick label style={font=\footnotesize},label style={font=\small},
 title style={font=\small},legend style={font=\scriptsize,draw=none},
 legend pos=north east,grid=major,grid style={gray!15}]
\nextgroupplot[title={$\widehat H_N$}]
\addplot+[color=blue!70!black,mark=*] coordinates {
(8,0.12197200)
(16,0.07912850)
(32,0.04047388)
(64,0.02507533)
(128,0.01643404)
};
\addlegendentry{$H=0.875$}
\addplot+[color=red!65!black,mark=triangle*] coordinates {
(8,0.12033059)
(16,0.05338812)
(32,0.02319463)
(64,0.01622471)
(128,0.00686936)
};
\addlegendentry{$H=0.750$}
\addplot+[color=green!45!black,mark=square*] coordinates {
(8,0.10165875)
(16,0.05576474)
(32,0.03028510)
(64,0.01554327)
(128,0.00828949)
};
\addlegendentry{$H=0.625$}
\nextgroupplot[title={$\widetilde\theta_N$},ylabel={}]
\addplot+[color=blue!70!black,mark=*] coordinates {
(8,0.30925080)
(16,0.24870846)
(32,0.18184441)
(64,0.12778967)
(128,0.09957850)
};
\addlegendentry{$H=0.875$}
\addplot+[color=red!65!black,mark=triangle*] coordinates {
(8,0.33336053)
(16,0.18331775)
(32,0.09979022)
(64,0.07794963)
(128,0.03834392)
};
\addlegendentry{$H=0.750$}
\addplot+[color=green!45!black,mark=square*] coordinates {
(8,0.26698091)
(16,0.18708602)
(32,0.12983445)
(64,0.07658462)
(128,0.04643009)
};
\addlegendentry{$H=0.625$}
\end{groupplot}
\end{tikzpicture}
\caption{RMSEs of the joint estimators as the observation resolution
increases, with $B=50$ and $n=512$.}
\label{fig:simulation-joint-rmse}
\end{figure}
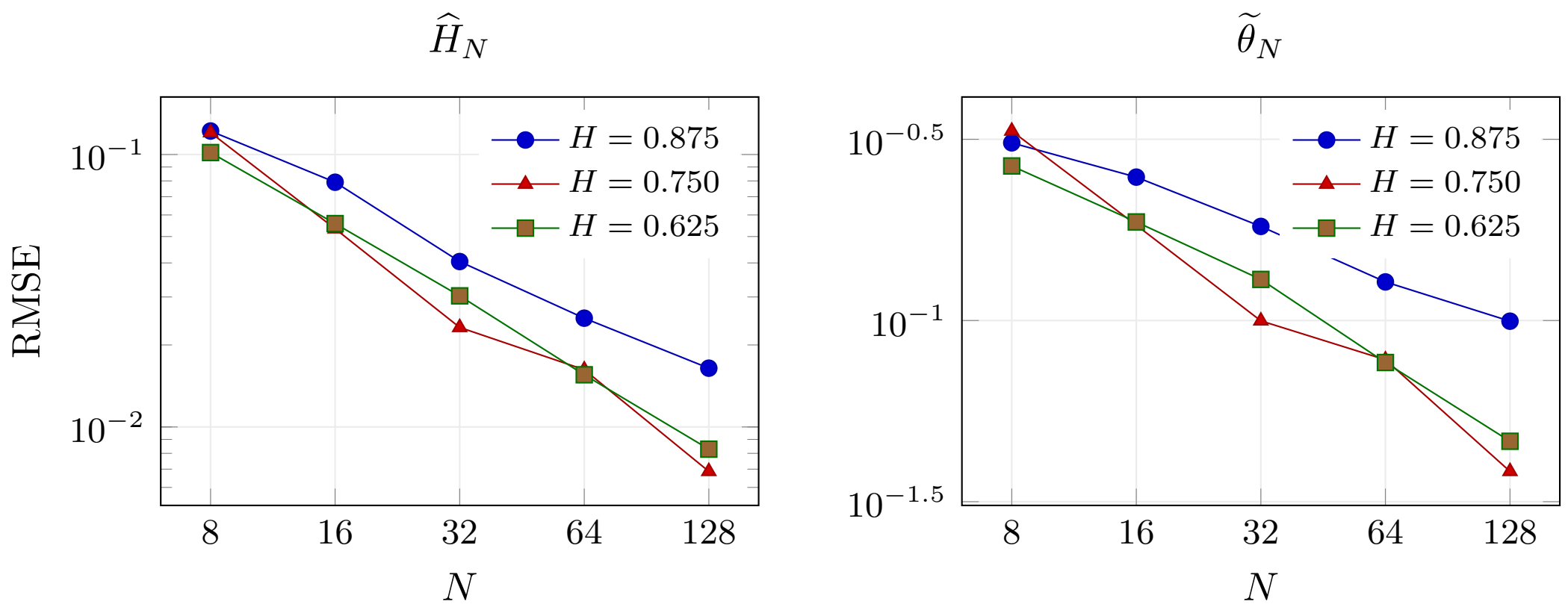

Figure~\ref{fig:simulation-joint-rmse} shows that the RMSEs of both
estimators decrease as $N$ increases. At $N=128$, negative empirical
biases remain visible for $H=0.625$. The projection is inactive in
all realizations with $N\ge32$.
Since the simulation resolution $n$ is fixed, the reported errors
include discretization effects. These results illustrate finite-sample
behavior; they do not establish sharp asymptotic rates for the
estimators.

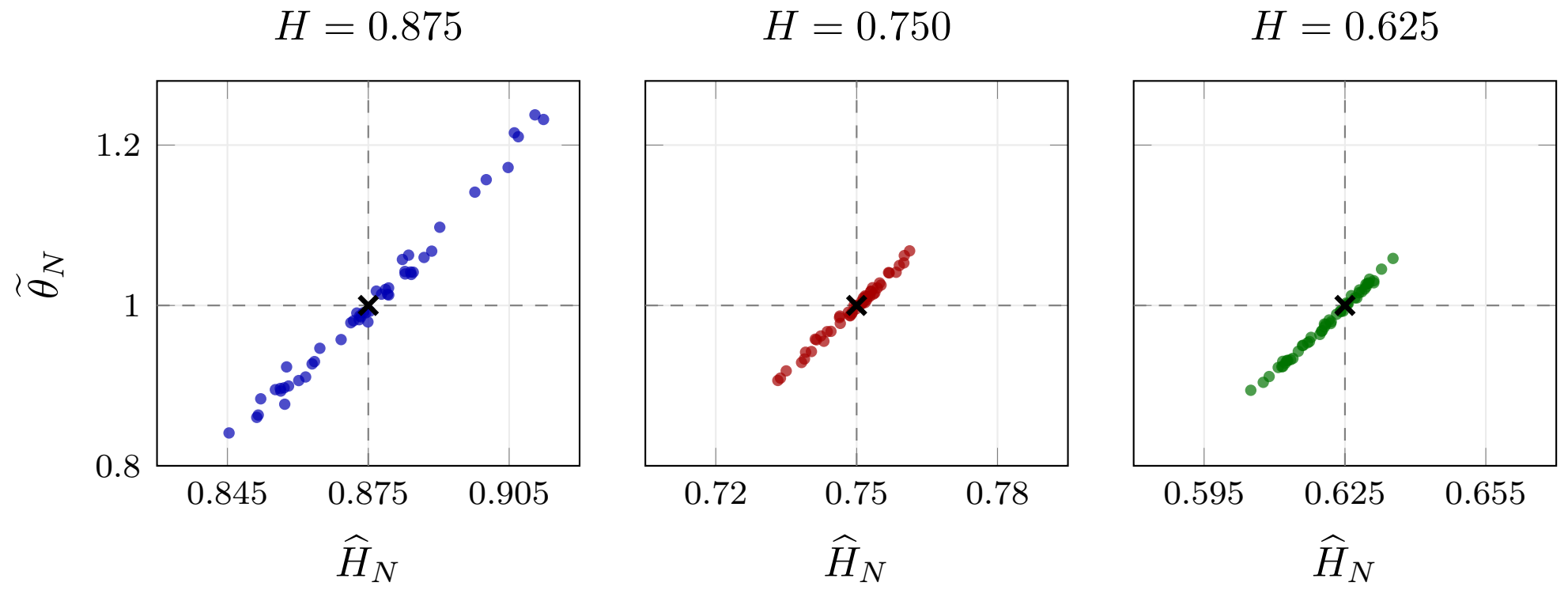
\begin{figure}[tb]
\centering
\begin{tikzpicture}
\begin{groupplot}[
 group style={group size=3 by 1,horizontal sep=0.55cm},
 width=0.32\textwidth,height=4.8cm,
 xlabel={$\widehat H_N$},ymin=0.8,ymax=1.28,
 ytick={0.8,1,1.2},
 tick label style={font=\scriptsize},label style={font=\small},
 title style={font=\small},
 scaled x ticks=false,xticklabel style={/pgf/number format/fixed,
 /pgf/number format/precision=3},grid=major,grid style={gray!15}]
\nextgroupplot[title={$H=0.875$},xmin=0.830,xmax=0.920,xtick={0.845,0.875,0.905},ylabel={$\widetilde\theta_N$}]
\addplot[gray,dashed] coordinates {(0.830,1) (0.920,1)};
\addplot[gray,dashed] coordinates {(0.875,0.8) (0.875,1.28)};
\addplot[only marks,mark=*,mark size=1.15pt,color=blue!70!black,opacity=0.7] coordinates {
(0.85625420,0.89642142)
(0.86164794,0.91078657)
(0.87255198,0.99052332)
(0.85707009,0.89724454)
(0.86467885,0.94656283)
(0.87445977,0.99077411)
(0.88358320,1.06272967)
(0.89022290,1.09751284)
(0.87305921,0.98221101)
(0.84534248,0.84095737)
(0.87532156,0.99279394)
(0.87394197,0.99095472)
(0.88687563,1.05985338)
(0.91231811,1.23194064)
(0.87934830,1.01297166)
(0.88226309,1.05730697)
(0.87184966,0.98089024)
(0.88280102,1.04237950)
(0.90011758,1.15686707)
(0.87910545,1.01395488)
(0.88415459,1.03866261)
(0.88400901,1.04166567)
(0.88284718,1.03923885)
(0.90477483,1.17203227)
(0.87671117,1.01787816)
(0.87932437,1.02185361)
(0.86018462,0.90632272)
(0.87315825,0.98689495)
(0.90694731,1.21040760)
(0.87127484,0.97839166)
(0.85759500,0.92332481)
(0.85636192,0.89324453)
(0.85156147,0.86332578)
(0.86353922,0.92998743)
(0.85124666,0.86038053)
(0.87341484,0.98576124)
(0.85801754,0.89951357)
(0.86920937,0.95740984)
(0.87491838,0.97941976)
(0.88850253,1.06772784)
(0.88457160,1.04144933)
(0.87781094,1.01403468)
(0.85209283,0.88355203)
(0.86303334,0.92705943)
(0.90610841,1.21529345)
(0.89769616,1.14128457)
(0.87864288,1.01990258)
(0.85720306,0.87679281)
(0.85520824,0.89500458)
(0.91049546,1.23770118)
};
\addplot[only marks,mark=x,mark size=3pt,very thick,black] coordinates {(0.875,1)};
\nextgroupplot[title={$H=0.750$},xmin=0.705,xmax=0.795,xtick={0.720,0.750,0.780},yticklabels={,,}]
\addplot[gray,dashed] coordinates {(0.705,1) (0.795,1)};
\addplot[gray,dashed] coordinates {(0.750,0.8) (0.750,1.28)};
\addplot[only marks,mark=*,mark size=1.15pt,color=red!65!black,opacity=0.7] coordinates {
(0.75360085,1.01418748)
(0.75279400,1.01130768)
(0.74642567,0.98678415)
(0.74374717,0.96734301)
(0.74858976,0.98772390)
(0.74453500,0.96767204)
(0.75213803,1.00697476)
(0.75694036,1.04054786)
(0.74145582,0.95712994)
(0.75489129,1.02799173)
(0.73322142,0.90645335)
(0.74037211,0.94251650)
(0.76015167,1.06220911)
(0.74984881,0.99562621)
(0.75178856,1.00359575)
(0.75440353,1.02173018)
(0.75203310,1.01111283)
(0.75234348,1.01061965)
(0.75519770,1.02537025)
(0.76128242,1.06806780)
(0.73888738,0.93320188)
(0.73500251,0.91837257)
(0.74823054,0.99147487)
(0.75143448,1.00416701)
(0.74892603,0.98944533)
(0.74302356,0.95531559)
(0.75387386,1.01562460)
(0.74991091,0.99931531)
(0.74651888,0.97763464)
(0.76005723,1.05303395)
(0.75142077,1.00928986)
(0.74944708,1.00020368)
(0.75108861,1.00508534)
(0.74635254,0.98493671)
(0.73826914,0.92883975)
(0.75841222,1.04138463)
(0.75112267,1.00221202)
(0.75182399,1.01198686)
(0.75909772,1.04979886)
(0.73912115,0.94164877)
(0.75077070,1.00090794)
(0.75307271,1.01742014)
(0.74118467,0.95792729)
(0.75679966,1.04097084)
(0.75338963,1.02210404)
(0.74236889,0.96191765)
(0.74861986,0.98706826)
(0.75300022,1.01758061)
(0.74902924,0.99104528)
(0.73377486,0.90938002)
};
\addplot[only marks,mark=x,mark size=3pt,very thick,black] coordinates {(0.750,1)};
\nextgroupplot[title={$H=0.625$},xmin=0.580,xmax=0.670,xtick={0.595,0.625,0.655},yticklabels={,,}]
\addplot[gray,dashed] coordinates {(0.580,1) (0.670,1)};
\addplot[gray,dashed] coordinates {(0.625,0.8) (0.625,1.28)};
\addplot[only marks,mark=*,mark size=1.15pt,color=green!45!black,opacity=0.7] coordinates {
(0.62753555,1.00941784)
(0.61263439,0.93132846)
(0.62918838,1.02038112)
(0.61593698,0.94998060)
(0.60494553,0.89418884)
(0.61743880,0.95491577)
(0.61177285,0.92382900)
(0.60885487,0.91133370)
(0.62857680,1.01700586)
(0.61174864,0.93025589)
(0.62091717,0.97619945)
(0.63010256,1.02812749)
(0.62712928,1.00947062)
(0.61078290,0.92260151)
(0.62406206,0.99330523)
(0.62063338,0.97313378)
(0.63523479,1.05862857)
(0.62799100,1.01587896)
(0.62163420,0.98184740)
(0.61704664,0.95346005)
(0.62321459,0.98900196)
(0.63278668,1.04527548)
(0.62550890,1.00141065)
(0.60760118,0.90409366)
(0.63114033,1.03081290)
(0.62980434,1.02723244)
(0.62197671,0.97741757)
(0.61391840,0.93396492)
(0.61156213,0.92360348)
(0.62463601,0.99298987)
(0.61341855,0.93235307)
(0.61232474,0.92799817)
(0.62571767,1.00245954)
(0.61283440,0.93127433)
(0.62412053,0.99228061)
(0.62953328,1.02329370)
(0.62810037,1.01932787)
(0.61965710,0.96384627)
(0.62945714,1.02640425)
(0.62207233,0.97976360)
(0.62526631,1.00292021)
(0.63023570,1.03278162)
(0.61999765,0.96786892)
(0.61778096,0.96017134)
(0.61506125,0.94276838)
(0.62062215,0.97685381)
(0.62636032,1.01221760)
(0.62015293,0.96839863)
(0.63114029,1.02823441)
(0.61612713,0.95024769)
};
\addplot[only marks,mark=x,mark size=3pt,very thick,black] coordinates {(0.625,1)};
\end{groupplot}
\end{tikzpicture}
\caption{Joint estimates from $50$ realizations at $N=128$ and
$n=512$. The crosses mark the true parameter pairs.}
\label{fig:simulation-joint}
\end{figure}

Figure~\ref{fig:simulation-joint} displays the joint estimates at
$N=128$. With truncation inactive, the identity in the proof of
Corollary~\ref{cor:joint-consistency} gives the first-order approximation
\[
 \widetilde\theta_N-\theta
 \approx (\widehat\theta_N-\theta)
 +\theta\left(\log N-\tfrac12\log2+\frac{1}{2H+1}\right)
       (\widehat H_N-H),
\]
where $\widehat\theta_N$ uses the true $H$.
At $\theta=1$ and $N=128$, the second term is approximately
$4.9(\widehat H_N-H)$; a Hurst estimation error of $0.01$ thus
contributes about $0.049$ to the diffusion estimation error.
This term dominates the fluctuations in these simulations,
explaining the near-linear pattern and strong positive correlation.

\FloatBarrier
\section*{Statements and Declarations}
\noindent\textbf{Funding.}
The authors received no funding for this work.

\smallskip
\noindent\textbf{Competing interests.}
The authors declare that they have no competing interests.

\smallskip
\noindent\textbf{Data availability.}
The simulation code and data supporting the numerical results
are available from the corresponding author upon request.

\begingroup
\renewcommand{\bibnumfmt}[1]{#1.}

\endgroup

\end{document}